\documentclass[11pt]{amsart}
\usepackage{amscd,amssymb,amsmath}
\usepackage[arrow,matrix]{xy}
\usepackage{graphicx}
\usepackage{cite}
\usepackage{ifpdf}
\usepackage{tikz}

\newtheorem{thm}[subsection]{Theorem}

\newtheorem{proposition}[subsection]{Proposition}
\newtheorem{cor}[subsection]{Corollary}

\newtheorem{rk}[subsection]{Remark}
\newtheorem{defn}[subsection]{Definition}

\numberwithin{equation}{section} 

\newcommand{\M}{{\mathcal M}}
\newcommand{\A}{{\mathcal A}}

\newcommand{\bea}{\begin{eqnarray}}
\newcommand{\eea}{\end{eqnarray}}

\newcommand{\Z}{\mathbb{Z}}

\newcommand{\R}{\mathbb{R}}

\begin{document}
\title [ A chain of evolution algebras]
{ A chain of evolution algebras}

\author {J.M. Casas, M. Ladra, U.A. Rozikov}

\address{J.\ M.\ Casas\\ Department of Applied Mathematics, E.U.I.T. Forestal, Pontevedra, University of Vigo, 36005, Spain.}
 \email {jmcasas@uvigo.es}
 \address{M.\ Ladra\\ Department of Algebra, University of Santiago de Compostela, 15782, Spain.}
 \email {manuel.ladra@usc.es}
 \address{U.\ A.\ Rozikov\\ Institute of mathematics and information technologies,
Tashkent, Uzbekistan.} \email {rozikovu@yandex.ru}

\begin{abstract} We introduce a notion of chain of evolution algebras.
The sequence of matrices of the structural constants for this
chain of evolution algebras satisfies an analogue of
Chapman-Kolmogorov equation. We give several examples (time
homogenous, time non-homogenous, periodic, etc.) of such chains.
For a periodic chain of evolution algebras we construct a
continuum set of non-isomorphic evolution algebras and show that
the corresponding discrete time chain of evolution algebras is
dense in the set. We obtain a criteria for an evolution algebra to
be baric and give a concept of a property transition. For several
chains of evolution algebras we describe the behavior of the baric
property depending on the time. For a chain of evolution algebras
given by the matrix of a two-state evolution we define a baric
property controller function and under some conditions on this
controller we prove that the chain is not baric almost surely
(with respect to Lebesgue measure). We also construct examples of
the almost surely baric chains of evolution algebras. We show that
there are chains of evolution algebras such that if it has a
unique (resp. infinitely many) absolute nilpotent element at a
fixed time, then it has unique (resp. infinitely many) absolute
nilpotent element any time; also there are chains of evolution
algebras which have not such property. For an example of two
dimensional chain of evolution algebras we give the full set of
idempotent elements and show that for some values of parameters
the number of idempotent elements does not depend on time, but for
other values of parameters there is a critical time $t_{c}$ such
that the chain has only two idempotent elements if time $t\geq
t_{\rm c}$ and it has four
idempotent elements if time $t< t_{\rm c}$.\\[3mm]

{\it AMS classifications (2010):} 17D92; 17D99; 60J27\\[2mm]

\end{abstract}
\maketitle
\section{Introduction} \label{sec:intro}

In this paper we consider some classes of non-associative algebras.
There exist several classes of non-associative algebras (baric,
evolution, Bernstein, train, stochastic, etc.), whose investigation
has provided a number of significant contributions to theoretical
population genetics. Such classes have been defined at different
times by several authors, and all the algebras belonging to these
classes are generally called ``genetic''. Etherington introduced the
formal language of abstract algebra to study of genetics in his
series of seminal papers \cite{e1,e2,e3}. In recent years many
authors have tried to investigate the difficult problem of
classification of these algebras. The most comprehensive references
for the mathematical research done in this area are \cite{ly,m,t,w}.

 In \cite{ly} an evolution algebra $\mathcal A$ associated to the
free population is introduced and using this non-associative algebra
many results are obtained in explicit form, e.g. the explicit
description of stationary quadratic operators, and the explicit
solutions of a nonlinear evolutionary equation in the absence of
selection, as well as general theorems on convergence to equilibrium
in the presence of selection.

In \cite{t} a new type of evolution algebra is introduced. This
evolution algebra is defined as follows.
 Let $(E,\cdot)$ be an algebra over a field $K$. If it admits a
basis $e_1,e_2,\dots$, such that $e_i\cdot e_j=0$, if $i\ne j$ and $
e_i\cdot e_i=\sum_{k}a_{ik}e_k$, for any $i$, then this algebra is
called an {\it evolution algebra}.

In this paper by the term {\it evolution algebra} we will understand
a finite dimensional evolution algebra $E$ (as mentioned above) over
field $\R$.

Evolution algebras have the following elementary properties (see
\cite{t}): Evolution algebras are not associative, in general; they
are commutative, flexible, but not power-associative, in general;
direct sums of evolution algebras are also evolution algebras;
Kronecker products of evolutions algebras are also evolution
algebras.

 The concept of evolution algebras lies between algebras and dynamical
systems. Algebraically, evolution algebras are non-associative
Banach algebra; dynamically, they represent discrete dynamical
systems. Evolution algebras have many connections with other
mathematical fields including graph theory, group theory, stochastic
processes, mathematical physics, etc.

In the book \cite{t}, the foundation of evolution algebra theory and
applications in non-Mendelian genetics and Markov chains are
developed, with pointers to some further research topics.

In \cite{rt} the algebraic structures of function spaces defined by
graphs and state spaces equipped with Gibbs measures by associating
evolution algebras are studied. Results of \cite{rt} also allow  a
natural introduction of thermodynamics in studying of several
systems of biology, physics and mathematics by theory of evolution
algebras.

The paper is organized as follows. In Section we give main
definitions related to a chain of evolution algebras. Therein we
give several examples (time homogenous, time non-homogenous,
periodic, etc.) of such chains. For a periodic chain of evolution
algebras we construct a continuum set of non-isomorphic evolution
algebras and show that the corresponding discrete time chain of
evolution algebras is dense in the set. In Section 3 we obtain a
criteria for an evolution algebra to be baric. The concept of a
property transition is introduced in Section 4. This section  also
contains several chains of evolution algebras for which we
describe the behavior of the baric property depending on the time.
For a chain of evolution algebras given by the matrix of a
two-state evolution we define a baric property controller function
and under some conditions on this controller we prove that the
chain is not baric almost surely (with respect to Lebesgue
measure). We also construct examples of the almost surely baric
chains of evolution algebras. We show that there are chains of
evolution algebras such that if it has a unique (resp. infinitely
many) absolute nilpotent element at a fixed time, then it has
unique (resp. infinitely many) absolute nilpotent element any
time; also there are chains of evolution algebras which have not
such property. In the last subsection for an example of two
dimensional chain of evolution algebras we give the full set of
idempotent elements and show that for some values of parameters
the number of idempotent elements does not depend on time, but for
other values of parameters there is a critical time $t_{c}$ such
that the chain has only two idempotent elements if time $t\geq
t_{\rm c}$ and it has four idempotent elements if time $t< t_{\rm
c}$.

\section{Definition and examples of CEA}
Consider a family $\left\{E^{[s,t]}:\ s,t \in \R,\ 0\leq s\leq t
\right\}$ of $n$-dimensional evolution algebras over the field $\R$,
with basis $e_1,\dots,e_n$ and multiplication table
\begin{equation}\label{1}
 e_ie_i = M_i^{[s,t]}=\sum_{j=1}^na_{ij}^{[s,t]}e_j, \ \ i=1,\dots,n; \ \
e_ie_j =0,\ \ i\ne j.\end{equation} Here parameters $s,t$ are
considered as time.

Denote by
$\M^{[s,t]}=\left(a_{ij}^{[s,t]}\right)_{i,j=1,\dots,n}$-the matrix
of structural constants.

\begin{defn}\label{d1} A family $\left\{E^{[s,t]}:\ s,t \in \R,\ 0\leq s\leq t
\right\}$ of $n$-dimensional evolution algebras over the field $\R$
is called a chain of evolution algebras (CEA) if the matrix
$\M^{[s,t]}$ of structural constants satisfies the
Chapman-Kolmogorov equation
\begin{equation}\label{2}
\M^{[s,t]}=\M^{[s,\tau]}\M^{[\tau,t]}, \ \ \mbox{for any} \ \
s<\tau<t.
\end{equation}
\end{defn}

If $\rho_i$ is a projection map of $E^{[s,t]}$, which maps every
element of $E^{[s,t]}$ to its $e_i$ component, then equation
(\ref{2}) can be written as
\begin{equation}\label{3}
M^{[s,t]}_i=\sum_{j=1}^n\rho_j(M^{[s,\tau]}_i)M^{[\tau,t]}_j, \ \
\mbox{for any} \ \ s<\tau<t.
\end{equation}
\begin{defn}\label{d2} A CEA is called a time-homogenous CEA if
the matrix $\M^{[s,t]}$ depends only on $t-s$. In this case we write
$\M^{[t-s]}$.
\end{defn}

\begin{defn}\label{dp} A CEA is called periodic if its matrix $\M^{[s,t]}$
is periodic with respect to at least one of the variables  $s$, $t$,
i.e. (periodicity with respect to $t$) $\M^{[s,t+P]}=\M^{[s,t]}$ for
all values of $t$. The constant $P$ is called the period, and is
required to be nonzero.
\end{defn}

\begin{rk}\label{r1}
 In general, an algebra $\A^{[s,t]}$ can be given by a cubic matrix
 $\M^{[s,t]}=\left(a_{ijk}^{[s,t]}\right)_{i,j,k=1,\dots,n}$ of structural constants.
 Our Definition \ref{d1} can be extended to $\A^{[s,t]}$ using analogues of the Chapman-Kolmogorov equations
 for quadratic operators (see \cite{gs,gp,gam}). Since in the general case there
 are two types of the Chapman-Kolmogorov equations: type $A$ and type $B$ \cite{gp}, one also can define two types of
 chain of (general) algebras using the Chapman-Kolmogorov equations of type $A$ and type
 $B$,
 respectively. In this paper we shall only consider CEA, which is more simple than general case, because it is defined
 by quadratic matrices.
\end{rk}

 {\sl The CEA corresponding to a Markov process.}

Let $\left\{\M^{[s,t]}, \ \ 0\leq s\leq t\right\}$ be a family of
stochastic matrices which satisfies the equation (\ref{2}), then it
defines a Markov process. Thus we have

\begin{thm}\label{t1} For each Markov process, there is a CEA whose
structural constants are transition probabilities of the process,
and whose generator set (basis) is the state space of the Markov
process.
\end{thm}

If $\M^{[s,t]}$ does not depend on time (i.e. $=\M$) then the CEA
contains only one evolution algebra $E$. Note that for a Markov
chain defined by $\M$ the corresponding $E$ has been studied in
\cite{t}.

Now we shall give several concrete examples of CEA.

{\bf Example} 1. To show a time dependent CEA we use the following
example of time homogenous Markov process (see \cite{k}) : for $n=3$
consider
$$a_{ii}^{[t]}={2\over 3}e^{-{3\over 2}At}\cos(\alpha t)+{1\over 3},
\ \ i=1,2,3;$$
$$a_{12}^{[t]}=a_{23}^{[t]}=a_{31}^{[t]}=e^{-{3\over 2}At}\left({1\over \sqrt{3}}\sin(\alpha t)-{1\over 3}\cos(\alpha t)\right)+{1\over 3};$$
$$a_{21}^{[t]}=a_{32}^{[t]}=a_{13}^{[t]}=-e^{-{3\over 2}At}\left({1\over \sqrt{3}}\sin(\alpha t)+{1\over 3}\cos(\alpha t)\right)+{1\over 3},$$
where $A>0$, $\alpha={\sqrt{3}\over 2}A$.

Let $E^{[t]}, t\geq 0$ be the corresponding CEA. It is easy to see
that $E^{[t]}$ has an oscillation behavior depending on time $t$.
Moreover $\lim_{t\to +\infty}E^{[t]}=E$, where $E$ is an evolution
algebra with the multiplication table
$$e_1^2=e_2^2=e_3^2={1\over 3}(e_1+e_2+e_3), \ \ e_ie_j=0, i\ne j.$$

{\sl The CEA corresponding to a family of matrices which do not
define a process.}

{\bf Example} 2. We shall give a time homogenous CEA which are
different from CEAs corresponding to Markov processes. For $n=2$
take $$a_{11}^{[t]}=a_{22}^{[t]}=a^{[t]}; \ \
a_{12}^{[t]}=a_{21}^{[t]}=b^{[t]}.$$ Then equation (\ref{2}) is
equivalent to
$$a^{[t]}=a^{[\tau]}a^{[t-\tau]}+b^{[\tau]}b^{[t-\tau]};$$
$$b^{[t]}=a^{[\tau]}b^{[t-\tau]}+b^{[\tau]}a^{[t-\tau]}.$$
Denote $f(t)=a^{[t]}+b^{[t]}, \ \ \varphi(t)=a^{[t]}-b^{[t]}$, then
the last system of functional equations can be written as
$$f(t)=f(\tau)f(t-\tau), \ \
\varphi(t)=\varphi(\tau)\varphi(t-\tau).$$ Both these equations are
known as exponential Cauchy equation and the system of equations has
solution $f(t)=\lambda^t$, $\varphi(t)=\mu^t$, where $\lambda,
\mu\geq 0$. Consequently,
 $a^{[t]}={1\over 2}(\lambda^t+\mu^t)$, $b^{[t]}={1\over 2}(\lambda^t-\mu^t)$. But this solution does not define any Markov
process, in general.

Let $E^{[t]}, t\geq 0$ be the corresponding CEA. Depending on
parameters $\lambda$ and $\mu$ we get distinct behavior of $E^{[t]}$
for $t\to +\infty$, i.e. we have
$$\lim_{t\to +\infty}E^{[t]}=\left\{\begin{array}{lll}
E_0 \ \ \mbox{if} \ \ 0<\lambda, \mu<1,\\[2mm]
E_1 \ \ \mbox{if} \ \ \lambda=\mu=1,\\[2mm]
E_{1/2} \ \ \mbox{if} \ \ \lambda=1, 0\leq\mu<1,\\[2mm]
E_{-1/2} \ \ \mbox{if} \ \ \mu=1, 0\leq\lambda<1,\\[2mm]
E_\infty \ \ \mbox{otherwise},\\[2mm]
\end{array}\right.
$$
where $E_0$ is an evolution algebra with zero multiplication; $E_1$
is an evolution algebra with multiplication table
$$e_1^2=e_1, e_2^2=e_2, \ \ e_1e_2=0;$$ $E_{1/2}$ is an evolution
algebra with multiplication table
$$e_1^2=e_2^2={1\over 2}(e_1+e_2), \ \ e_1e_2=0;$$
$E_{-1/2}$ is an evolution algebra with multiplication table
$$e_1^2={1\over 2}(e_1-e_2), \ \ e_2^2=-{1\over 2}(e_1-e_2),\ \ e_1e_2=0;$$
and
 $E_\infty$ is a vector space which has ``infinity
multiplication'', or we can say that in $E_\infty$ an algebra
structure is not defined. This example shows that a limit of a CEA
can be non evolution algebra.

{\bf Example} 3. {\it A two-state evolution.}  Now we shall give an
example of time non-homogeneous CEA, the matrix of structural
constants of which also does not define any (time non homogenous)
Markov process in general. Consider $n=2$ and matrix
$\M^{[s,t]}=\left(a_{ij}^{[s,t]}\right)_{i,j=1,2}$ with
\begin{equation}\label{aa}
\begin{array}{ll}
a_{11}^{[s,t]}={1\over 2}\left(1+\alpha(s,t)+\beta(s,t)\right),\ \
a_{12}^{[s,t]}={1\over
2}\left(1-\alpha(s,t)-\beta(s,t)\right),\\[3mm]
a_{21}^{[s,t]}={1\over 2}\left(1+\alpha(s,t)-\beta(s,t)\right),\ \
a_{22}^{[s,t]}={1\over
2}\left(1-\alpha(s,t)+\beta(s,t)\right).\\[3mm]
\end{array}
\end{equation}
 In
this case the equation (\ref{2}) is equivalent to (see \cite{ht})
\begin{equation}\label{4}
\begin{array}{ll}
\alpha(s,t)=\alpha(\tau,t)+\alpha(s,\tau)\beta(\tau,t),\\[2mm]
\beta(s,t)=\beta(s,\tau)\beta(\tau,t), \ \ s<\tau<t.\\[2mm]
\end{array}
\end{equation}
The second equation of the system (\ref{4}) is known as Cantor's
second equation, it has very rich family of solutions:
$\beta(s,t)={\Phi(t)\over \Phi(s)}$, where $\Phi$ is an arbitrary
function with $\Phi(s)\ne 0$. Using this function $\beta$ for the
function $\alpha$ we obtain
$$
{\alpha(s,t)\over \Phi(t)} ={\alpha(\tau,t)\over \Phi(t)}
+{\alpha(s,\tau)\over \Phi(\tau)}.$$ Now denote
$\gamma(s,t)={\alpha(s,t)\over \Phi(t)}$ then the last equation gets
the following form
$$\gamma(s,t)=\gamma(s,\tau)+\gamma(\tau,t).$$
This equation is known as Cantor's first equation which also has
very rich family of solutions: $\gamma(s,t)=\Psi(t)-\Psi(s)$, where
$\Psi$ is an arbitrary function. Hence a solution
$\M^{[s,t]}=\left(a_{ij}^{[s,t]}\right)_{i,j=1,2}$ to the equation
(\ref{2}) is given by
$$a_{11}^{[s,t]}={1\over 2}\left(1+\Phi(t)(\Psi(t)-\Psi(s))+{\Phi(t)\over \Phi(s)}\right),$$
$$a_{12}^{[s,t]}={1\over 2}\left(1-\Phi(t)(\Psi(t)-\Psi(s))-{\Phi(t)\over \Phi(s)}\right),$$
$$a_{21}^{[s,t]}={1\over 2}\left(1+\Phi(t)(\Psi(t)-\Psi(s))-{\Phi(t)\over \Phi(s)}\right),$$
$$a_{22}^{[s,t]}={1\over 2}\left(1-\Phi(t)(\Psi(t)-\Psi(s))+{\Phi(t)\over \Phi(s)}\right).$$
Let $E^{[s,t]}$, $0\leq s\leq t$ be the corresponding to this
solution CEA. This CEA varies by two parameters, for example, if
$t=s$ we get $E^{[t,t]}=E$ with multiplication table $e^2_1=e_1,
e_2^2=e_2, e_1e_2=0$. Moreover, choosing functions $\Phi$ and $\Psi$
one can variate the limit behavior of the CEA. For example, if
$\Phi$ and $\Psi$ such that $\lim_{t\to +\infty}\Phi(t)\Psi(t)=
\lim_{t\to +\infty}\Phi(t)=0$, then for a fixed $s$ we have
$\lim_{t\to +\infty}E^{[s,t]}=E_{1/2}$, where $E_{1/2}$ is an
evolution algebra with multiplication table
$$e_1^2=e_2^2={1\over 2}(e_1+e_2), \ \ e_1e_2=0.$$

{\bf Example} 4. {\it A $n$-dimensional time non-homogenous CEA.}
Here for arbitrary $n$ we shall give an example of time
non-homogenous CEA. Let $\{A^{[t]}, \,t\geq 0\}$ be a family of
invertible (for all $t$), $n\times n$ matrices. Define the following
matrix
$$M^{[s,t]}=A^{[s]}(A^{[t]})^{-1},$$
where $(A^{[t]})^{-1}$ is the inverse of $A^{[t]}$.

This matrix satisfies the equation (\ref{2}). Indeed, using
associativity of the multiplication of matrices we get
$$\M^{[s,\tau]}\M^{[\tau,t]}=A^{[s]}
\left((A^{[\tau]})^{-1}A^{[\tau]}\right)(A^{[t]})^{-1}=A^{[s]}(A^{[t]})^{-1}=\M^{[s,t]}.$$
Thus each family (with one parameter) of invertible $n\times n$
matrices defines a CEA $E^{[s,t]}$ which is time non-homogenous, in
general. But will be a time homogenous CEA, for example, if
$A^{[t]}$ is equal to $t$th power of an invertible matrix $A$.

Construction of a family of invertible $n\times n$ matrices
$A^{[t]}$ is not difficult, for example, one can take $A^{[t]}$ as a
triangular $n\times n$ matrix of the form
$$
A^{[t]}  =\left(\begin{array}{ccccc}
a_{11}^{[t]}  & 0            &0       & \dots          & 0 \\[3mm]
a_{21}^{[t]}  & a_{22}^{[t]} &    0    &  \dots               &0 \\[3mm]
a_{31}^{[t]}  & a_{32}^{[t]} & \ddots &   \dots              &0\\[3mm]
\vdots        & \vdots       & \ddots &\ddots          &0 \\[3mm]
a_{n1}^{[t]}  & a_{n2}^{[t]} & \ldots & a_{nn-1}^{[t]} &
a_{nn}^{[t]}
\end{array}\right).
$$
which is called lower triangular matrix or one can take an upper
triangular matrix. Then the matrices are invertible iff
$a_{ii}^{[t]}\ne 0$, for all $i=1,\dots,n$ and $t$. So this example
also gives a very rich class of CEAs.

  {\bf Example} 5. {\it Periodic CEA.} To get a
periodic CEA, we can consider the $E^{[s,t]}$ constructed in Example
3, and choose $\Phi$ and $\Psi$ as periodic (non-constant)
functions. Then corresponding CEA is periodic. In this case  for any
fixed $s$, the limit $\lim_{t\to+\infty}E^{[s,t]}$ does not exist in
general, moreover its set of limit points (evolution algebras) can
be a continuum set. We shall make this point clear as follows.
Construct a time homogenous CEA which is periodic. Consider $n=2$
take
$$a_{11}^{[t]}=a_{22}^{[t]}=a^{[t]}; \ \
a_{12}^{[t]}=-b^{[t]}, \ \ a_{21}^{[t]}=b^{[t]}.$$ Then equation
(\ref{2}) is equivalent to
$$a^{[t]}=a^{[\tau]}a^{[t-\tau]}-b^{[\tau]}b^{[t-\tau]};$$
$$b^{[t]}=a^{[\tau]}b^{[t-\tau]}+b^{[\tau]}a^{[t-\tau]}.$$
This system reminds the following identities
$$\cos t=\cos \tau \cos(t-\tau)-\sin \tau \sin(t-\tau);$$
$$\sin t=\cos \tau \sin(t-\tau)+\sin \tau \cos(t-\tau).$$
Consequently, one solution
$\M^{[t]}=\left(a_{ij}^{[t]}\right)_{i,j=1,2}$ to equation (\ref{2})
is
\begin{equation}\label{m}
\M^{[t]}=\left(\begin{array}{cc}
\cos t &  \sin t\\
-\sin t &  \cos t
\end{array}\right).
\end{equation}
Since that matrix is periodic with period $P=2\pi$, the
corresponding CEA $E^{[t]}$ is also periodic. Moreover this CEA is
very interesting: for arbitrary 2-dimensional evolution algebra
$E^+_a$, or $E^-_a$, $a\in [-1,1]$ with structural constants
matrix
$$\M^{\pm}_a=\left(\begin{array}{cc}
a &  \pm\sqrt{1-a^2}\\
\mp\sqrt{1-a^2} & a
\end{array}\right)$$
respectively, there is a sequence $t_n=t_n(a)$ of times such that
$\lim_{n\to\infty}E^{[t_n]}=E^+_a$ or $E^-_a$. We have
$E^{\pm}_a\ne E^{\pm}_b$ if $a\ne b$. Moreover the following is
true

\begin{proposition}\label{pp} 1) For any $a,b\in [-1,1], \, a\ne \pm b$,
the algebras $E^+_a$ and $E^+_b$ are not isomorphic. The algebras
$E^+_{a}$ and $E^+_{-a}$ are isomorphic.

2) For any $a,b\in [-1,1], \, a\ne \pm b$, the algebras $E^-_a$
and $E^-_b$ are not isomorphic. The algebras $E^-_{a}$ and
$E^-_{-a}$ are isomorphic.
\end{proposition}

\proof 1) Let $\varphi=\left(\begin{array}{cc}
\alpha &  \beta\\[3mm]
\gamma & \delta\\[3mm]
\end{array}\right)$ be an isomorphism of
the evolution algebra $E^+_a$ to the evolution algebra $E^+_b$.
Here $\det(\varphi)\ne 0$. By the multiplication table of the
evolution algebras, we get the following relation between matrices
$\M^+_a$ and $\M^+_b$:
$$\M^+_b={1\over \det(\varphi)}\times$$ $$\left(\begin{array}{cc}
(a\delta-\sqrt{1-a^2}\gamma)\alpha^2-(a\gamma+\sqrt{1-a^2}\delta)\beta^2&
(a\alpha+\sqrt{1-a^2}\beta)\beta^2-(a\beta-\sqrt{1-a^2}\alpha)\alpha^2
\\[3mm]
(a\delta-\sqrt{1-a^2}\gamma)\gamma^2-(a\gamma+\sqrt{1-a^2}\delta)\delta^2&
(a\alpha+\sqrt{1-a^2}\beta)\delta^2-(a\beta-\sqrt{1-a^2}\alpha)\gamma^2
\\[3mm]
\end{array}\right).$$

Since $\det(\M^+_a)=1$, it is easy to see that there are two
classes of isomorphisms:
$$\mathcal C_1=\left\{\left(\begin{array}{cc}
\alpha & 0\\
0& \delta
\end{array}\right): \alpha\delta\ne 0\right\}, \ \
\mathcal C_2=\left\{\left(\begin{array}{cc}
0 & \beta \\
\gamma & 0
\end{array}\right): \beta\gamma\ne 0\right\}.$$

For the class $\mathcal C_1$ the matrix $\M^+_b$ must satisfy the
following
$$\M^+_b=\left(\begin{array}{cc}
b &  \sqrt{1-b^2}\\
-\sqrt{1-b^2} & b
\end{array}\right)=\left(\begin{array}{cc}
a\alpha &  \sqrt{1-a^2}{\alpha^2\over \delta}\\[3mm]
-\sqrt{1-a^2}{\delta^2\over \alpha} & a\delta
\end{array}\right).$$
From this equality we get $\alpha=\delta=\sqrt{{1-b^2\over
1-a^2}}={b\over a}$ if $a\ne 0,\pm 1$ which is satisfied iff
$a=b$. Hence the isomorphisms from the class $\mathcal C_1$ can
not give an isomorphism from $\M^+_a$ to $\M^+_b$. For $a=0$ we
get $b=0$. One can take $\alpha=\delta=\mp 1$ if $a=\pm 1$ and
$b=\mp 1$. Hence $E^+_{\pm 1}$ is isomorph to $E^+_{\mp 1}$.

For the class $\mathcal C_2$ the matrix $\M^+_b$ must satisfy the
following
$$\M^+_b=\left(\begin{array}{cc}
b &  \sqrt{1-b^2}\\
-\sqrt{1-b^2} & b
\end{array}\right)=\left(\begin{array}{cc}
a\beta &  -\sqrt{1-a^2}{\beta^2\over \gamma}\\[3mm]
\sqrt{1-a^2}{\gamma^2\over \beta} & a\gamma
\end{array}\right).$$
From this equality we get $\beta=\gamma=-\sqrt{{1-b^2\over
1-a^2}}={b\over a}$ if $a\ne 0, \pm 1$ which is satisfied iff $a=-
b$. Hence the isomorphisms from the class $\mathcal C_2$ can only
give an isomorphism from $\M^+_a$ to $\M^+_{-a}$.

2) The proof of 2) is similar to the proof of 1).
\endproof
Consider now discrete time $n$, $n\in \mathbb N$ and the CEA
$\{E^{[n]}, n\in \mathbb N\}$ given by matrix (\ref{m}).

\begin{proposition}\label{pp} The discrete time CEA $E^{[n]}$, $n\in \mathbb N$, is dense in
the set $\{E^{\pm}_a, \, a\in [-1,1]\}$ of evolution algebras,
i.e. for an arbitrary evolution algebra $E^{\pm}_a$ there exists a
sequence $\{n_k\}_{k=1,2,...}$ of natural numbers such that
$\lim_{k\to\infty}E^{[n_k]}=E^+_a$ or  $E^-_a$.
\end{proposition}
\proof It is known that the sequences $\{\sin n\}$ and $\{\cos
n\}, n\in \mathbb N$, are dense in $[-1,1]$ (see e.g.\cite{g}).
Hence for any $a\in [-1,1]$ there is a sequence
$\{n_k\}_{k=1,2,...}$ of natural numbers such that
$\lim_{k\to\infty}\cos(n_k)=a$. The same sequence can be used to
get $\lim_{k\to\infty}E^{[n_k]}=E^+_a$ or $E^-_a$.
\endproof

\section{A criterion for an evolution algebra to be baric}

A {\it character} for an algebra $A$ is a nonzero multiplicative
linear form on $A$, that is, a nonzero algebra homomorphism from $A$
to $\R$ \cite{ly}. Not every algebra admits a character. For
example, an algebra with the zero multiplication has no character.

\begin{defn}\label{d3} A pair $(A, \sigma)$ consisting of an algebra $A$ and a
character $\sigma$ on $A$ is called a {\it baric algebra}. The
homomorphism $\sigma$ is called the weight (or baric) function of
$A$ and $\sigma(x)$ the weight (baric value) of $x$.
\end{defn}
In \cite{ly} for the evolution algebra of a free population it is
proven that there is a character $\sigma(x)=\sum_i x_i$, therefore
that algebra is baric. But the evolution algebra $E$ introduced in
\cite{t} is not baric, in general. The following theorem gives a
criterion for an evolution algebra $E$ to be baric.

\begin{thm}\label{t2} An $n$-dimensional evolution algebra $E$, over field the $\R$,
is baric if and only if there is a column
$\left(a_{1i_0},\dots,a_{ni_0}\right)^T$ of its structural constants
matrix $\M=\left(a_{ij}\right)_{i,j=1,\dots,n}$, such that
$a_{i_0i_0}\ne 0$ and $a_{ii_0}=0$, for all $i\ne i_0$. Moreover,
the corresponding weight function is $\sigma(x)=a_{i_0i_0}x_{i_0}$.
\end{thm}
\proof {\sl Necessity.} Take $x,y\in E$ with $x=\sum_{i=1}^nx_ie_i$,
$y=\sum_{i=1}^ny_ie_i$. Assume $\sigma(x)=\sum_{i=1}^n \alpha_i
x_i$, $ x\in E$ is a character. We have
$$\sigma(xy)=\sum_{i=1}^n\left(\sum_{j=1}^na_{ij}\alpha_j\right)x_iy_i;
\ \ \sigma(x)\sigma(y)=\sum_{i=1}^n\sum_{j=1}^n
\alpha_i\alpha_jx_iy_j.$$ From $\sigma(xy)=\sigma(x)\sigma(y)$ we
get
\begin{equation}\label{5}
\alpha_i\alpha_j=0\ \ \mbox{for any } \ \ i\ne j, \ i,j=1,\dots,n;
\end{equation}
\begin{equation}\label{6}
\sum_{j=1}^n a_{ij}\alpha_j=\alpha^2_i\ \ \mbox{for any } \ \
i=1,\dots,n.
\end{equation}

It is easy to see that the system (\ref{5}) has a solution
$\alpha=(\alpha_1,\dots,\alpha_n)$ with
$\alpha_1^2+\dots+\alpha_n^2> 0$ if and only if exactly one
coordinate of $\alpha$, say $\alpha_{i_0}$, is not zero, and all
others are zeros. Substituting this solution in (\ref{6}) we get
$$
\begin{array}{ll}
a_{ii_0}\alpha_{i_0}=0, \ \ \mbox{if} \ \ i\ne i_0, \ i=1,\dots,n;\\[2mm]
a_{i_0i_0}\alpha_{i_0}=\alpha^2_{i_0}, \ \ \mbox{if} \ \ i=i_0.
\end{array}
$$
From the last equations we get $a_{i_0i_0}\ne 0$, $a_{ii_0}=0$, for
all $i\ne i_0$ and $\alpha_{i_0}=a_{i_0i_0}$.

{\sl Sufficiency.} Assume there is a column
$\left(a_{1i_0},\dots,a_{ni_0}\right)^T$, such that $a_{i_0i_0}\ne
0$ and $a_{ii_0}=0$, for all $i\ne i_0$. Then it is easy to see that
$\sigma(x)=a_{i_0i_0}x_{i_0}$ is a weight function, therefore $E$ is
a baric evolution algebra.
\endproof

A baric algebra $A$ may have several weight functions. As a
corollary of Theorem \ref{t2} we have

\begin{cor} If the matrix $\M$, mentioned in Theorem \ref{t2}, has
several columns $\left(a_{1i_j},\dots,a_{ni_j}\right)^T$, \ \
$j=i_1,\dots,i_m$, $m\leq n$, which satisfy conditions of Theorem
\ref{t2} then the evolution algebra $E$ has exactly $m$ weight
functions $\sigma(x)=a_{i_ji_j}x_{i_j}$, $j=i_1,\dots,i_m$.
\end{cor}

There are two types of trivial evolution algebras \cite{t}: {\it
zero evolution algebra}, which satisfies $e_ie_j=0$ for all
$i,j=1,\dots,n$; {\it non-zero trivial evolution algebra}, which
satisfies $e_ie_j=0$ for all $i\ne j$ and $e_i^2=a_{ii}e_i$, where
$a_{ii}\in \R$ is non-zero for some $i=1,\dots,n$. By Theorem
\ref{t2} we conclude that the zero evolution algebra is not baric,
but any non-zero trivial evolution algebra is a baric algebra.
Moreover, there are baric evolution algebras which are not trivial.

\section{Property transition}

If a system has parameters (as usually like: temperature, time,
interaction, etc.) then a property of the system can variate by a
parameter. For example, the behavior of phases (states) of a system
in physics, depends on temperature $T>0$, if for some values of $T$
there is a unique phase and for other values there are several
phases, then the physical system has a phase transition \cite{ge}.
Similar transitions of a property can be seen for systems of
biology, chemistry, etc. Here we shall define a notion of property
transition for CEA.

\begin{defn}\label{d4} Assume a CEA, $E^{[s,t]}$, has a property, say $P$,
at pair of times $(s_0,t_0)$; we say that the CEA has $P$ property
transition if there is a pair $(s,t)\ne (s_0,t_0)$ at which the CEA
has no the property $P$.
\end{defn}

Denote
$$\mathcal T=\{(s,t): 0\leq s\leq t\};$$
$$\mathcal T_P=\{(s,t)\in \mathcal T: E^{[s,t]} \ \ \mbox{has property} \  P \};$$
$$\mathcal T_P^0=\mathcal T\setminus \mathcal T_P=\{(s,t)\in \mathcal T: E^{[s,t]} \ \ \mbox{has no property} \ P \}.$$

\begin{defn}\label{d5} We call the set

$\mathcal T_P$-the duration of the property $P$;

$\mathcal T_P^0$-the lost duration of the property $P$;

The partition $\{\mathcal T_P, \mathcal T^0_P\}$ of the set
$\mathcal T$ is called $P$ property diagram.
\end{defn}

For example, if $P=$commutativity then since any evolution algebra
is commutative, we conclude that any CEA has not commutativity
property transition.

\subsection{Baric property transition.}

Since a CEA is not a baric algebra, in general, using Theorem
\ref{t2} we can give baric property diagram. Let us do this for the
above given Examples 1-4.

 {\bf Example} 1'. For the case of Example
1, by Theorem \ref{t2} we have that $E^{[t]}$ is baric iff
$$a_{ii}^{[t]}={2\over
3}e^{-{3\over 2}At}\cos(\alpha t)+{1\over 3}=1.$$ This has unique
solution $t=0$. Consequently, $\mathcal T_{\rm baric}=\{0\}$,
$\mathcal T^0_{\rm baric}=\{t: t>0\}$. Thus the CEA $E^{[t]}$ is
baric (even non zero trivial) evolution algebra only at initial
time, and it loses baricity as soon as the time turned on.

{\bf Example} 2'. In Example 2, using Theorem \ref{t2} we obtain
that

$$\mathcal T_{\rm baric}=\begin{cases}
\{0\} & \mbox{if} \ \ \lambda\ne \mu;\\[2mm]
\mathcal T &\mbox{if} \ \ \lambda=\mu.
\end{cases}
$$
Thus the CEA $E^{[t]}$ has not baric property transition if
$\lambda=\mu$, and it has a baric property transition, as in Example
1', if $\lambda\ne \mu$.

{\bf Example} 3'. {\it Baric property transition for a two-state
evolution.} Since in case of Example 3, we have a rich class of CEA
here we shall give a special theory of the baric property
transition. Using Theorem \ref{t2} we obtain that ${\mathcal T}_{\rm
baric}$ is the set of $(s,t)$ such that
$$1+\Phi(t)(\Psi(t)-\Psi(s))-{\Phi(t)\over \Phi(s)}=0\ \ \mbox{or}
\ \ 1-\Phi(t)(\Psi(t)-\Psi(s))-{\Phi(t)\over \Phi(s)}=0.$$
These equations can be rewritten as
$$ \theta(t)=\theta(s), \ \ \theta^-(t)=\theta^-(s),$$
where
\begin{equation}\label{7}
\theta(t)={1\over \Phi(t)}+\Psi(t), \, \theta^-(t)={1\over
\Phi(t)}-\Psi(t).
\end{equation}

Thus
$$\mathcal T_{\rm baric}=\mathcal T_{\rm baric}(\theta)\cup \mathcal T_{\rm baric}(\theta^-),$$
here $\mathcal T_{\rm baric}(\theta)= \left\{(s,t)\in \mathcal T:
\theta(t)=\theta(s)\right\}$.

\begin{rk} To describe the set $\mathcal T_{\rm baric}$ one has to describe the sets
$\mathcal T_{\rm baric}(\theta)$ and $\mathcal T_{\rm
baric}(\theta^-)$, both of which are defined by the parameter
functions $\Phi$ and $\Psi$. Note that if we replace $\Phi$ with
$-\Phi$ or $\Psi$ with $-\Psi$ then these sets transfer to each
other. Since $\Phi$ and $\Psi$ are arbitrary functions, it will be
enough to describe only $\mathcal T_{\rm baric}(\theta)$ for
arbitrary $\theta$. Thus in the sequel of this subsection we shall
deal with description of $\mathcal T_{\rm baric}(\theta)$.
\end{rk}

 The function $\theta(t)$ is called {\it baric property
controller} of the CEA. Because, it really controls the baric
duration set, for example, if $\theta$ is a strong monotone function
then the duration is ``minimal'', i.e. the line $s=t$, but if
$\theta$ is a constant function then the baric duration set is
``maximal'', i.e. it is $\mathcal T$. Since $\Phi$ and $\Psi$ are
arbitrary functions, we have a rich class of controller functions,
therefore we have a ``powerful''  control on the property to be
baric.

\begin{figure}[h]
\begin{center}
\scalebox{.8}{

\begin{tikzpicture}[x=1.00mm, y=1.00mm, inner xsep=0pt, inner ysep=0pt, outer xsep=0pt, outer ysep=0pt]
\path[line width=0mm] (73.34,31.94) rectangle +(71.92,69.73);
\definecolor{L}{rgb}{0,0,0}
\path[line width=0.30mm, draw=L] (139.38,40.38) -- (80.06,40.25);
\definecolor{F}{rgb}{0,0,0}
\path[line width=0.30mm, draw=L, fill=F] (139.38,40.38) --
(136.58,41.07) -- (139.38,40.38) -- (136.59,39.67) -- (139.38,40.38)
-- cycle; \path[line width=0.30mm, draw=L] (80.06,40.38) --
(80.44,40.38); \path[line width=0.60mm, draw=L] (129.82,89.88) --
(80.06,40.12); \path[line width=0.30mm, draw=L] (80.06,99.07) --
(80.06,40.00); \path[line width=0.30mm, draw=L, fill=F]
(80.06,99.07) -- (79.36,96.27) -- (80.06,99.07) -- (80.76,96.27) --
(80.06,99.07) -- cycle; \draw(137.60,35.55) node[anchor=base
west]{\fontsize{14.23}{17.07}\selectfont $t$}; \draw(75.34,95.75)
node[anchor=base west]{\fontsize{14.23}{17.07}\selectfont $s$};
\path[line width=0.60mm, draw=L] (99.83,59.90) .. controls
(102.54,54.66) and (103.12,52.18) .. (108.13,47.14) .. controls
(112.06,43.18) and (117.32,41.68) .. (121.65,40.51); \path[line
width=0.15mm, draw=L, dash pattern=on 2.00mm off 1.00mm]
(99.96,59.45) -- (99.96,40.12); \draw(98.50,35.55) node[anchor=base
west]{\fontsize{14.23}{17.07}\selectfont $t_{\rm c}$};
\draw(119.50,35.55) node[anchor=base
west]{\fontsize{14.23}{17.07}\selectfont $t'_{\rm c}$};
\draw(87.10,75.32) node[anchor=base
west]{\fontsize{14.23}{17.07}\selectfont Baric curves}; \path[line
width=0.15mm, draw=L] (125.33,50.13) -- (125.24,62.98); \path[line
width=0.15mm, draw=L, fill=F] (125.33,50.13) -- (126.01,52.93) --
(125.33,50.13) -- (124.61,52.92) -- (125.33,50.13) -- cycle;
\path[line width=0.15mm, draw=L] (108.93,48.06) -- (97.20,74.11);
\path[line width=0.15mm, draw=L, fill=F] (108.93,48.06) --
(108.42,50.90) -- (108.93,48.06) -- (107.14,50.32) -- (108.93,48.06)
-- cycle; \draw(116.18,64.19) node[anchor=base
west]{\fontsize{14.23}{17.07}\selectfont Non-baric set}; \path[line
width=0.15mm, draw=L] (103.41,65.66) -- (97.20,74.11); \path[line
width=0.15mm, draw=L, fill=F] (103.41,65.66) -- (102.32,68.33) --
(103.41,65.66) -- (101.19,67.50) -- (103.41,65.66) -- cycle;
\path[line width=0.15mm, draw=L] (107.03,42.10) -- (125.24,63.07);
\path[line width=0.15mm, draw=L, fill=F] (107.03,42.10) --
(109.40,43.76) -- (107.03,42.10) -- (108.34,44.68) -- (107.03,42.10)
-- cycle;
\end{tikzpicture}}
\end{center}
\caption[baric]{ The baric property diagram for
$\theta(t)=\lambda^{-t}+ct$.}\label{fig:1}
\end{figure}
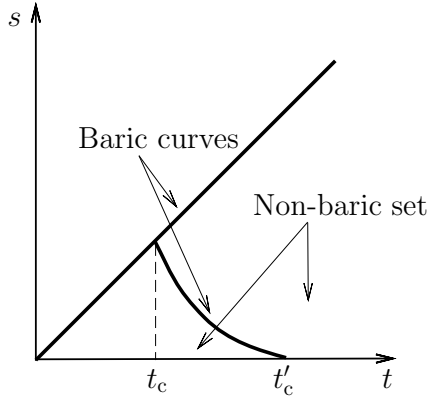

For a special choose of $\theta$ we have
\begin{proposition}\label{p1} If $\Phi(t)=\lambda^t$, $\lambda>0$ and $\Psi(t)=ct$, $c\in \R$. Then
$$
{\mathcal T}_{\rm baric}(\theta)={\mathcal T}_{\rm
baric}(\lambda,c)=\{(s,t): s=t\}\cup$$ $$
\begin{cases}
\ \ \ \ \emptyset & {\rm if} \ \
0<\lambda\leq 1, \, c\geq \ln\lambda; \ \ {\rm or} \\
& \lambda>1, \, c\in (-\infty,0]\cup[\ln\lambda,+\infty),\\
\{(s,t): 0\leq s\leq t_{\rm c}, \, t_c\leq t\leq t'_{\rm c}, \,
\theta(s)=\theta(t)\} & {\rm if} \ \
0<\lambda\leq 1, \, c<\ln\lambda; \ \ {\rm or}\\
& \lambda>1, \, c\in (0,\ln\lambda),\\
\end{cases}
$$
where $t_{\rm c}$ and $t'_{\rm c}$ serve as critical times, which
defined by $t_{\rm c}={1\over \ln\lambda}\ln\left({\ln\lambda\over
c}\right)$ and $t'_{\rm c}>0$ is a unique solution to
$\theta(t'_{\rm c})=1$.
\end{proposition}
\proof Under the conditions of the proposition we have
$\theta(t)=\lambda^{-t}+ct$, and the simple analysis of the equation
$\theta(s)=\theta(t)$ for this $\theta$ gives the full set
${\mathcal T}_{\rm baric}(\lambda,c)$.
\endproof
In Figure \ref{fig:1}, the baric property diagram is given

As a corollary of Proposition \ref{p1} we have

\begin{cor} 1) For any fixed $s$, with $0\leq s<t_c$ (resp. $t_c\leq s\leq t$), the time
$t$ has two (resp. one) critical values: $t_c^{(1)}=s$ (resp. $s$)
and $t_c^{(2)}$ which is a unique solution of
$\theta(t_c^{(2)})=\theta(s)$.

2) For any fixed $t$, with $0\leq s\leq t\leq t_c$ or $t'_c<t$
(resp. $t_c<t\leq t'_c$), the time $s$ has one (resp. two) critical
values: $s_c^{(1)}=t$ (resp. $s_c^{(1)}=t$ and $s_c^{(2)}$ which is
a unique solution of $\theta(t_c^{(2)})=\theta(s))$ .
\end{cor}

Let us discuss some more examples of the controller $\theta$. If
$\theta(t)=\tan(t)$ then $\tan(s)=\tan(t)$ has solution $t=s+\pi k$,
$k\in \Z$. The intersection of this family of lines with $\mathcal
T$ gives the family of half lines, i.e.
$$\mathcal T_{baric}(\tan(t))=\bigcup_{k=0,1,2,\dots}\left\{(s,t)\in \mathcal T: s=t-\pi
k\right\}.$$ If $\theta(t)=\sin(t)$ then $\sin(s)=\sin(t)$ has two
family of solutions: $s=t+2\pi k$, $k\in \Z$ and $s=-t+(2k+1)\pi$,
$k\in \Z$. The intersection of these families of lines with
$\mathcal T$ is
$$\mathcal T_{baric}(\sin(t))=\bigcup_{k=0,1,2,\dots}\left\{(s,t)\in \mathcal T: t=s+2\pi
k \ \ \mbox{or} \ \ t=-s+(2k+1)\pi\right\}.$$

In all above considered examples we obtained a set $\mathcal
T_{baric}(\theta)$ which has zero Lebesgue measure. But there is
controllers for which this set has non-zero Lebesgue measure, for
example, if  $\theta(t)$ is a controller function with the graph as
shown in Figure \ref{fig:2}, then the corresponding baric property
diagram is as shown in Figure \ref{fig:3}. Thus any ``constant
part'' of the graph of the controller gives a full triangle in the
diagram, moreover, any ``non-constant part'' gives several curves.
In this case the set $\mathcal T_{baric}(\theta)$ has a non-zero
Lebesgue measure.

\begin{figure}[h]
\begin{center}
\scalebox{.7}{

\begin{tikzpicture}[x=1.00mm, y=1.00mm, inner xsep=0pt, inner ysep=0pt, outer xsep=0pt, outer ysep=0pt]
\path[line width=0mm] (73.85,6.71) rectangle +(109.00,96.15);
\definecolor{L}{rgb}{0,0,0}
\path[line width=0.30mm, draw=L] (80.06,100.85) -- (79.93,19.84);
\definecolor{F}{rgb}{0,0,0}
\path[line width=0.30mm, draw=L, fill=F] (80.06,100.85) --
(79.35,98.06) -- (80.06,100.85) -- (80.75,98.05) -- (80.06,100.85)
-- cycle; \path[line width=0.30mm, draw=L] (79.93,20.09) --
(169.62,20.09); \path[line width=0.30mm, draw=L, fill=F]
(169.62,20.09) -- (166.82,20.79) -- (169.62,20.09) -- (166.82,19.39)
-- (169.62,20.09) -- cycle;
\path[line width=0.30mm, draw=L] (79.92,50.55) -- (130.,50.55); 
\path[line width=0.30mm, draw=L] (80.19,44.72) -- (133.,44.72); 
\path[line width=0.30mm, draw=L] (80.06,66.15) -- (138.5,66.15); 
\path[line width=0.30mm, draw=L] (80.06,84.6) -- (142.5,84.6); 
\path[line width=0.60mm, draw=L] (80.06,66.15) .. controls
(85.42,58.11) and (88.70,52.27) .. (97.66,51.10) .. controls
(106.03,50.00) and (112.11,61.35) .. (115.02,70.49) .. controls
(116.80,76.10) and (117.22,79.47) .. (119.74,83.25) .. controls
(122.77,87.80) and (125.84,78.88) .. (126.12,76.36) .. controls
(126.88,69.47) and (128.03,62.20) .. (129.05,54.67) .. controls
(129.26,53.11) and (131.87,42.77) .. (133.26,45.36) .. controls
(137.85,53.90) and (136.98,64.28) .. (139.90,75.47) .. controls
(140.66,78.40) and (141.43,81.33) .. (142.57,84.27); \path[line
width=0.60mm, draw=L] (142.57,84.27) .. controls (142.70,86.31) and
(144.09,86.72) .. (144.69,87.65); \path[line width=0.60mm, draw=L]
(144.55,87.65) -- (153.04,87.65); \path[line width=0.60mm, draw=L]
(153.04,87.65) .. controls (155.04,92.30) and (158.39,95.50) ..
(160.05,100.34); \draw(75.85,94.60) node[anchor=base
west]{\fontsize{14.23}{17.07}\selectfont $s$}; \draw(169.62,15.12)
node[anchor=base west]{\fontsize{14.23}{17.07}\selectfont $t$};
\path[line width=0.30mm, draw=L, dash pattern=on 2.00mm off 1.00mm]
(99.83,51.10) -- (99.96,19.97); \draw(98,15.50) node[anchor=base
west]{\fontsize{12}{15}\selectfont $t_1$}; \path[line width=0.30mm,
draw=L, dash pattern=on 2.00mm off 1.00mm] (113.61,66.28) --
(113.74,20.22); \draw(112.0,15.50) node[anchor=base
west]{\fontsize{12}{15}\selectfont $t_2$}; \path[line width=0.30mm,
draw=L, dash pattern=on 2.00mm off 1.00mm] (121.52,84.65) --
(121.65,20.22); \draw(120.5,15.50) node[anchor=base
west]{\fontsize{12}{15}\selectfont $t_3$}; \path[line width=0.30mm,
draw=L, dash pattern=on 2.00mm off 1.00mm] (132.75,44.97) --
(132.75,19.97); \path[line width=0.30mm, draw=L, dash pattern=on
2.00mm off 1.00mm] (127.01,66.41) -- (127.01,20.22);
\draw(124.5,15.50) node[anchor=base
west]{\fontsize{12}{15}\selectfont $t_4$}; \path[line width=0.30mm,
draw=L] (180.85,53.55); \path[line width=0.30mm, draw=L, dash
pattern=on 2.00mm off 1.00mm] (129.82,50.99) -- (129.82,20.12);
\draw(128.29,15.50) node[anchor=base
west]{\fontsize{12}{15}\selectfont $t_5$}; \draw(131.88,15.50)
node[anchor=base west]{\fontsize{12}{15}\selectfont $t_6$};
\path[line width=0.30mm, draw=L, dash pattern=on 2.00mm off 1.00mm]
(135.68,50.87) -- (135.68,19.86); \draw(135.4,15.50)
node[anchor=base west]{\fontsize{12}{15}\selectfont $t_7$};
\path[line width=0.30mm, draw=L, dash pattern=on 2.00mm off 1.00mm]
(138.36,66.18) -- (138.36,20.12); \draw(138.53,15.50)
node[anchor=base west]{\fontsize{12}{15}\selectfont $t_8$};
\path[line width=0.30mm, draw=L, dash pattern=on 2.00mm off 1.00mm]
(142.60,84.68) -- (142.60,20.12); \draw(142.,15.50) node[anchor=base
west]{\fontsize{12}{15}\selectfont $t_9$}; \path[line width=0.30mm,
draw=L, dash pattern=on 2.00mm off 1.00mm] (144.80,87.61) --
(144.80,19.86); \draw(145.3,15.50) node[anchor=base
west]{\fontsize{12}{15}\selectfont $t_{10}$}; \path[line
width=0.30mm, draw=L, dash pattern=on 2.00mm off 1.00mm]
(153.0,87.61) -- (153.0,19.86); \draw(152.5,15.50) node[anchor=base
west]{\fontsize{12}{15}\selectfont $t_{11}$};
\end{tikzpicture}
}
\end{center}
\caption{An example of controller $\theta$.} \label{fig:2}
\end{figure}
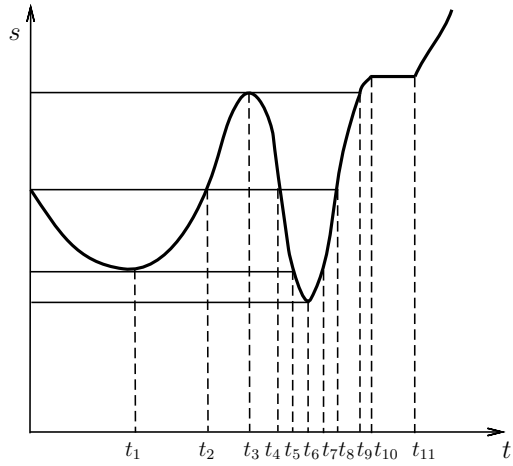

\begin{figure}[h]
\begin{center}
\scalebox{.7}{

\begin{tikzpicture}[x=1.00mm, y=1.00mm, inner xsep=0pt, inner ysep=0pt, outer xsep=0pt, outer ysep=0pt]
\path[line width=0mm] (52.72,-3.01) rectangle +(128.97,126.84);
\definecolor{L}{rgb}{0,0,0}
\path[line width=0.30mm, draw=L] (69.85,120.09) -- (69.85,9.88);
\definecolor{F}{rgb}{0,0,0}
\path[line width=0.30mm, draw=L, fill=F] (69.85,120.09) --
(69.15,117.29) -- (69.85,120.09) -- (70.55,117.29) -- (69.85,120.09)
-- cycle; \path[line width=0.30mm, draw=L] (69.57,10.29) --
(179.69,10.03); \path[line width=0.30mm, draw=L, fill=F]
(179.69,10.03) -- (176.89,10.74) -- (179.69,10.03) -- (176.89,9.34)
-- (179.69,10.03) -- cycle; \path[line width=0.30mm, draw=L]
(70.80,10.15) -- (70.53,10.15); \path[line width=0.60mm, draw=L]
(69.57,10.15) -- (169.88,110.04); \path[line width=0.60mm, draw=L]
(102.69,42.75) .. controls (105.14,40.57) and (108.01,34.89) ..
(110.19,33.); \path[line width=0.60mm, draw=L] (110.19,33.0) ..
controls (113.94,29.41) and (117.52,25.15) .. (119.40,21.0);
\path[line width=0.60mm, draw=L] (119.40,21.0) .. controls
(118.23,15.36) and (112.91,13.17) .. (110.19,10.20); \path[line
width=0.60mm, draw=L] (124.22,64.27) .. controls (126.10,62.51) and
(128.45,61.68) .. (130.10,59.64); \path[line width=0.60mm, draw=L]
(130.10,21.5) .. controls (132.94,16.61) and (139.38,12.47) ..
(143.10,10.20); \path[line width=0.60mm, draw=L] (143.10,32.94) ..
controls (136.76,29.94) and (134.85,26.53) .. (131.04,22.31) ..
controls (130.81,22.06) and (129.71,21.36) .. (130.10,21.5);
\path[line width=0.60mm, draw=L] (130.10,42.74) .. controls
(134.17,39.83) and (139.25,34.23) .. (142.80,32.90); \path[line
width=0.60mm, draw=L] (142.80,50.64) .. controls (145.31,47.13) and
(150.42,45.57) .. (152.8,43.0);
\definecolor{F}{rgb}{0.502,0.502,0.502}
\path[fill=F] (169.74,109.80) -- (169.88,99.86) -- (159.93,99.72) --
(169.74,109.67) -- (169.61,109.40) -- cycle;
\definecolor{L}{rgb}{0.502,0.502,0.502}
\path[line width=0.15mm, draw=L]  (160.00,99.72) -- (160.00,99.79)
(162.00,99.75) -- (162.00,101.82) (164.00,99.77) -- (164.00,103.85)
(166.00,99.80) -- (166.00,105.87) (168.00,99.83) -- (168.00,107.90);
\definecolor{L}{rgb}{0,0,0}
\path[line width=0.60mm, draw=L] (169.74,109.80) -- (169.88,99.86)
-- (159.93,99.72) -- (169.74,109.67) -- (169.61,109.40);
\draw(175.47,5.82) node[anchor=base
west]{\fontsize{14.23}{17.07}\selectfont $t$}; \draw(78,5.5)
node[anchor=base west]{\fontsize{12}{15}\selectfont $t_1$};
\draw(90.,5.5) node[anchor=base west]{\fontsize{12}{15}\selectfont
$t_2$}; \draw(100.5,5.5) node[anchor=base
west]{\fontsize{12}{15}\selectfont $t_3$}; \draw(108.5,5.5)
node[anchor=base west]{\fontsize{12}{15}\selectfont $t_4$};
\draw(117.,5.5) node[anchor=base west]{\fontsize{12}{15}\selectfont
$t_5$}; \draw(122.,5.5) node[anchor=base
west]{\fontsize{12}{15}\selectfont $t_6$}; \draw(128.,5.5)
node[anchor=base west]{\fontsize{12}{15}\selectfont $t_7$};
\draw(142.5,5.5) node[anchor=base west]{\fontsize{12}{15}\selectfont
$t_8$}; \draw(151.,5.5) node[anchor=base
west]{\fontsize{12}{15}\selectfont $t_9$}; \draw(160.,5.5)
node[anchor=base west]{\fontsize{12}{15}\selectfont $t_{10}$};
\draw(168.10,5.5) node[anchor=base
west]{\fontsize{12}{15}\selectfont $t_{11}$}; \draw(60,19)
node[anchor=base west]{\fontsize{12}{15}\selectfont $t_1$};
\draw(60,32.5) node[anchor=base west]{\fontsize{12}{15}\selectfont
$t_2$}; \draw(60,41.) node[anchor=base
west]{\fontsize{12}{15}\selectfont $t_3$}; \draw(60,48.5)
node[anchor=base west]{\fontsize{12}{15}\selectfont $t_4$};
\draw(60,59) node[anchor=base west]{\fontsize{12}{15}\selectfont
$t_5$}; \draw(60,63.5) node[anchor=base
west]{\fontsize{12}{15}\selectfont $t_6$}; \draw(60,69.60)
node[anchor=base west]{\fontsize{12}{15}\selectfont $t_7$};
\draw(60,82.) node[anchor=base west]{\fontsize{12}{15}\selectfont
$t_8$}; \draw(60,90.5) node[anchor=base
west]{\fontsize{12}{15}\selectfont $t_9$}; \draw(60,99)
node[anchor=base west]{\fontsize{12}{15}\selectfont $t_{10}$};
\draw(60,109) node[anchor=base west]{\fontsize{12}{15}\selectfont
$t_{11}$}; \draw(64.95,117.91) node[anchor=base
west]{\fontsize{12}{15}\selectfont {$s$}};
\definecolor{L}{rgb}{0.502,0.502,0.502}
\path[line width=0.15mm, draw=L] ;
\definecolor{L}{rgb}{0,0,0}
\path[line width=0.30mm, draw=L, dash pattern=on 2.00mm off 1.00mm]
(91.9,32.93) -- (91.9,10);
\definecolor{L}{rgb}{0.502,0.502,0.502}
\path[line width=0.15mm, draw=L] ;
\definecolor{L}{rgb}{0,0,0}
\path[line width=0.30mm, draw=L, dash pattern=on 2.00mm off 1.00mm]
(102.42,42.33) -- (102.42,10.);
\definecolor{L}{rgb}{0.502,0.502,0.502}
\path[line width=0.15mm, draw=L] ;
\definecolor{L}{rgb}{0.502,0.502,0.502}
\path[line width=0.15mm, draw=L] ;
\definecolor{L}{rgb}{0,0,0}
\path[line width=0.30mm, draw=L, dash pattern=on 2.00mm off 1.00mm]
(110.3,50.64) -- (110.3,10.);
\definecolor{L}{rgb}{0.502,0.502,0.502}
\path[line width=0.15mm, draw=L] ;
\definecolor{L}{rgb}{0,0,0}
\path[line width=0.30mm, draw=L, dash pattern=on 2.00mm off 1.00mm]
(119.45,60.18) -- (119.45,10);
\definecolor{L}{rgb}{0.502,0.502,0.502}
\path[line width=0.15mm, draw=L] ;
\definecolor{L}{rgb}{0,0,0}
\path[line width=0.30mm, draw=L, dash pattern=on 2.00mm off 1.00mm]
(123.9,64.2) -- (123.9,9.9);
\definecolor{L}{rgb}{0,0,0}
\path[line width=0.30mm, draw=L, dash pattern=on 2.00mm off 1.00mm]
(129.8,69.59) -- (129.8,10);
\definecolor{L}{rgb}{0.502,0.502,0.502}
\path[line width=0.15mm, draw=L] ;
\definecolor{L}{rgb}{0,0,0}
\path[line width=0.30mm, draw=L, dash pattern=on 2.00mm off 1.00mm]
(142.8,83.21) -- (142.8,10);
\definecolor{L}{rgb}{0.502,0.502,0.502}
\path[line width=0.15mm, draw=L] ;
\definecolor{L}{rgb}{0,0,0}
\path[line width=0.30mm, draw=L, dash pattern=on 2.00mm off 1.00mm]
(152.50,92.62) -- (152.50,10);
\definecolor{L}{rgb}{0.502,0.502,0.502}
\path[line width=0.15mm, draw=L] ;
\definecolor{L}{rgb}{0,0,0}
\path[line width=0.30mm, draw=L, dash pattern=on 2.00mm off 1.00mm]
(160,99.57) -- (160,10.0);
\definecolor{L}{rgb}{0.502,0.502,0.502}
\path[line width=0.15mm, draw=L] ;
\definecolor{L}{rgb}{0,0,0}
\path[line width=0.30mm, draw=L, dash pattern=on 2.00mm off 1.00mm]
(169.7,99.81) -- (169.7,10.);
\definecolor{L}{rgb}{0.502,0.502,0.502}
\path[line width=0.15mm, draw=L] ;
\definecolor{L}{rgb}{0,0,0}
\path[line width=0.30mm, draw=L, dash pattern=on 2.00mm off 1.00mm]
(69.85,21) -- (152.50,21.);
\definecolor{L}{rgb}{0.502,0.502,0.502}
\path[line width=0.15mm, draw=L] ;
\definecolor{L}{rgb}{0,0,0}
\path[line width=0.30mm, draw=L, dash pattern=on 2.00mm off 1.00mm]
(80.3,20.81) -- (80.30,9.77);
\definecolor{L}{rgb}{0.502,0.502,0.502}
\path[line width=0.15mm, draw=L] ;
\definecolor{L}{rgb}{0,0,0}
\path[line width=0.30mm, draw=L, dash pattern=on 2.00mm off 1.00mm]
(69.85,43.) -- (152.50,43.0);
\definecolor{L}{rgb}{0.502,0.502,0.502}
\path[line width=0.15mm, draw=L] ;
\definecolor{L}{rgb}{0,0,0}
\path[line width=0.30mm, draw=L, dash pattern=on 2.00mm off 1.00mm]
(70.53,32.8) -- (152.50,32.80);
\definecolor{L}{rgb}{0.502,0.502,0.502}
\path[line width=0.15mm, draw=L] ;
\definecolor{L}{rgb}{0,0,0}
\path[line width=0.30mm, draw=L, dash pattern=on 2.00mm off 1.00mm]
(69.71,50.25) -- (152.50,50.25);
\definecolor{L}{rgb}{0.502,0.502,0.502}
\path[line width=0.15mm, draw=L] ;
\definecolor{L}{rgb}{0,0,0}
\path[line width=0.30mm, draw=L, dash pattern=on 2.00mm off 1.00mm]
(69.98,100.) -- (159.79,100);
\definecolor{L}{rgb}{0.502,0.502,0.502}
\path[line width=0.15mm, draw=L] ;
\definecolor{L}{rgb}{0,0,0}
\path[line width=0.30mm, draw=L, dash pattern=on 2.00mm off 1.00mm]
(69.71,109.99) -- (169.74,109.99); \path[line width=0.60mm, draw=L]
(80.34,20.80) .. controls (83.44,16.08) and (87.02,12.62) ..
(91.85,10.3); \path[line width=0.60mm, draw=L] (152.8,43.) ..
controls (151.40,40.58) and (147.65,38.01) .. (146.30,36.47) ..
controls (145.35,35.38) and (142.8,34.53) .. (142.8,33.06);
\path[line width=0.30mm, draw=L, dash pattern=on 2.00mm off 1.00mm]
(69.57,60.05) -- (152.50,60.18); \path[line width=0.30mm, draw=L,
dash pattern=on 2.00mm off 1.00mm] (69.85,64.41) -- (152.50,64.00);
\path[line width=0.30mm, draw=L, dash pattern=on 2.00mm off 1.00mm]
(69.85,91.79) -- (152.50,92.06); \path[line width=0.30mm, draw=L,
dash pattern=on 2.00mm off 1.00mm] (69.85,83.34) -- (152.50,82.79);
\path[line width=0.30mm, draw=L, dash pattern=on 2.00mm off 1.00mm]
(69.57,70.00) -- (152.50,70.00);
\end{tikzpicture}
}
\end{center}
\caption{The baric property diagram for the controller $\theta$ with
graph as in Figure \ref{fig:2}.} \label{fig:3}
\end{figure}
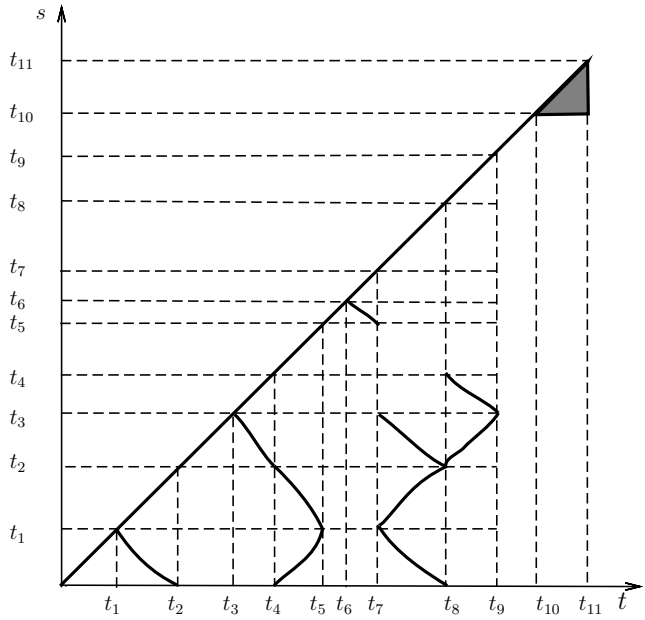

Let $\theta(t)=D(t)$ be the Dirichlet function defined by
$$D(t)=\left\{\begin{array}{ll}
1 \ \ \mbox{if} \ \  t \, \mbox{rational};\\
0 \ \ \mbox{if} \ \  t \, \mbox{irrational}.\\
\end{array}\right.$$
In this case we have a rich set of baric property duration, i.e.
$$\mathcal T_{baric}(D(t))=\{(s,t)\in \mathcal T: t \, {\rm and} \, s \, \mbox{rational}\}\cup
\{(s,t)\in \mathcal T: t \, {\rm and} \, s \,\mbox{irrational}\}.$$

\begin{defn}\label{d6} A function $\theta$ defined on $\R$ is called a
function of {\it countable variation} if it has the following
properties:

1. it is continuous except at most on a countable set, (which is
denoted by $X_c=\{x_1,x_2,\dots\}$), it has only jump-type
discontinuities (denote the one-sided limit from the negative
direction by $\theta(x_i^-)$ and from the positive direction by
$\theta(x_i^+)$, $i=1,2,\dots$);

2.  it has at most a countable set of singular (extremum) points
(which is denoted by $X_e=\{y_1,y_2,\dots\}$).
\end{defn}

Note that any function of countable variation has not ``constant
parts'' in its graph.

The following theorem gives a characteristics of the baric property
duration set.

\begin{thm}\label{t3} If the controller $\theta$ (see (\ref{7})) is a function of
countable variation,  then the baric duration set $\mathcal T_{\rm
baric}(\theta)$ has zero Lebesgue measure, that is the corresponding
CEA is not baric almost surely.
\end{thm}
\proof Using the (finite or infinite) sequences $X_c$ and $X_e$ we
construct the sequences $\{t^-_{i,k}\}_{k=1,2,\dots}$, $i=1,2,\dots$
with $\theta(t^-_{i,k})=\theta(x^-_i)$ for all $k$;
$\{t^+_{i,q}\}_{q=1,2,\dots}$, $i=1,2,\dots$ with
$\theta(t^+_{i,q})=\theta(x^+_i)$ and
$\{t^{e}_{j,l}\}_{l=1,2,\dots}$, $j=1,2,\dots$, where
$\theta(t^{e}_{j,l})=\theta(y_j)$ for all $l$. Now define a sequence
$\{t_i\}_{i=1,2,\dots}$, with $t_1<t_2<t_3<\dots$ as follows
$$\{t_i\}_{i=1,2,\dots}=X_c\cup X_e\bigcup_{i}\left(
\{t^-_{i,k}\}_{k=1,2,\dots}\cup
\{t^+_{i,q}\}_{q=1,2,\dots}\right)\bigcup_{j}\{t^{e}_{j,l}\}_{l=1,2,\dots}.$$
Since $\theta$ is a function of countable variation, the sequence
$\{t_i\}_{i=1,2,\dots}$ is at most countable. In a case, if it is a
bounded sequence (in particular, a finite sequence), then we add the
last term to be $+\infty$. Consider rectangles
$$\mathcal T_{ij}=\left\{(s,t)\in \R^2: t_i\leq s\leq t_{i+1}, \, t_j\leq t\leq
t_{j+1}\right\}.$$

Denote $G(\theta)=\{(t,y): y=\theta(t)\}.$

 By the construction, the rectangles have the following
properties:

-- The set of all rectangles is at most a countable set;

-- The intersection $G(\theta)\cap \mathcal T_{ij}$ is empty or
contains a monotone part of the graph $G(\theta)$.

If $G(\theta)\cap \mathcal T_{ij}$ is empty then we say $\mathcal
T_{ij}$ is empty.

Now we shall construct the set $\mathcal T_{\rm baric}(\theta)$. Fix
$i, j$ such that the rectangle $\mathcal T_{ij}$ is not-empty (an
empty rectangle does not give any contribution to the set $\mathcal
T_{\rm baric}(\theta)$), then we have
 $$\mathcal T_{\rm baric}(\theta)\cap\mathcal T_{kj}=\mbox{a curve
 giving an
 one-to-one corespondence between} \ \ [t_j,t_{j+1}]$$ $$ \mbox{and}
 \ \ [t_k,t_{k+1}]\ \ \mbox{if} \ \ \mathcal T_{ik}\ne \emptyset,
 \,k=1,\dots,j-1.$$
Thus we have
$$\mathcal T_{\rm baric}(\theta)=\bigcup_{kj}\left(\mathcal T_{\rm baric}(\theta)\cap\mathcal
T_{kj}\right).$$ Since there are a countable set of rectangles and
in each rectangle we may have at most a curve which has Lebesgue
measure zero (because, these curves give one-to-one
correspondences), we conclude that the set $\mathcal T_{\rm
baric}(\theta)$ also has zero Lebesgue measure.
\endproof

{\bf Example} 4'. Consider the CEA $E^{[s,t]}$ constructed in
Example 4 by a family of invertible lower (or upper) triangular
matrices $A^{[t]}$, $t\geq 0$.

\begin{thm}\label{t4}  For any pair of time $(s,t)$ the $n$-dimensional evolution algebra
$E^{[s,t]}$, constructed by a family of (lower or upper) triangular
invertible matrices is baric. Moreover, $E^{[s,t]}$ has a weight
function $\sigma(x)=\M^{[s,t]}_{nn}x_{n}$, where $\M^{[s,t]}_{ii}$,
$i=1,\dots,n$ are diagonal entries of
$\M^{[s,t]}=A^{[s]}(A^{[t]})^{-1}$.
\end{thm}
\proof It is known that the standard operations on triangular
matrices conveniently preserve the triangular form: the sum and
product of two lower triangular matrices is again lower triangular.
The inverse of a lower triangular matrix is also lower triangular,
and of course we can multiply a lower triangular matrix by a
constant and it will still be lower triangular. This means that the
lower triangular matrices form a subalgebra of the ring of square
matrices for any given size. The analogous result holds for upper
triangular matrices. Using these properties we get that $\M^{[s,t]}$
is also a triangular matrix. Moreover, since $A^{[t]}$ is
invertible, its determinant is non-zero for all $t$. Thus
$$\det(\M^{[s,t]})=\prod_{i=1}^n\M^{[s,t]}_{ii}=
\det(A^{[s]})\det((A^{[t]})^{-1})\ne 0.$$ Consequently, all diagonal
entries of the matrix are non-zero. In particular,
$\M^{[s,t]}_{nn}\ne 0$, and Theorem \ref{t2} completes the proof.
\endproof
\begin{cor} The CEA $E^{[s,t]}$ constructed by triangular
invertible matrices has not baric property transition. \end{cor}

\subsection{Absolute nilpotent elements transition.}

The element $x$ of an algebra $A$ is called an {\it absolute
nilpotent} if $x^2=0$.

Let $E=\R^n$ be an evolution algebra over the field $\R$ with
structural constant coefficients matrix $\M=(a_{ij})$, then for
arbitrary $x=\sum_ix_ie_i$ and $y=\sum_iy_ie_i\in \R^n$ we have
$$xy=\sum_j\left(\sum_ia_{ij}x_iy_i\right)e_j, \ \
x^2=\sum_j\left(\sum_ia_{ij}x^2_i\right)e_j.$$

For a $n$-dimensional evolution algebra $\R^n$ consider operator
$V \colon \R^n\to \R^n$, $x\mapsto V(x)=x'$ defined as
\begin{equation}\label{v}
x'_j= \sum_{i=1}^na_{ij}x_i^2, \ \ j=1,\dots,n.
\end{equation}
This operator is called {\it evolution operator} \cite{ly}.

We have $V(x)=x^2$, hence the equation $V(x)=x^2=0$ is given by the
following system
\begin{equation}\label{n1}
\sum_ia_{ij}x_i^2=0, \ \ j=1,\dots,n.
\end{equation}

If $\det(\M)\ne 0$ then the system (\ref{n1}) has unique solution
$(0,\dots,0)$. If  $\det(\M)=0$ and rank$(\M)=r$ then we can
assume that the first $r$ rows of $\M$ are linearly independent,
consequently, the system of equations (\ref{n1}) can be written as
\begin{equation}\label{n11}
x_i^2=-\sum_{j=r+1}^nd_{ij}x_j^2, \ \ i=1,\dots,r,
\end{equation}
where $d_{ij}={\det(\M_{ij})\over\det(M_r)}$ with
$\M_r=\left(a_{ij}\right)_{i,j=1,\dots,r}$,
$$
\M_{ij} =\left(\begin{array}{ccccccc}
a_{11} &\dots & a_{i-1,1} & a_{j1} & a_{i+1,1} & \dots & a_{r1} \\[3mm]
a_{12} &\dots & a_{i-1,2} & a_{j2} & a_{i+1,2} & \dots & a_{r2} \\[3mm]
 &\dots &  & \dots &  & \dots & \\[3mm]
a_{1r} &\dots & a_{i-1,r} & a_{jr} & a_{i+1,r} & \dots & a_{rr}
\end{array}\right).
$$
An interesting problem is to find a necessary and sufficient
condition on matrix $D=(d_{ij})_{{i=1,\dots,r\atop
j=r+1,\dots,n}}$ under which the system (\ref{n11}) has unique
solution. The difficulty of the problem depends on rank $r$, here
we shall consider the case $r=n-1$.

\begin{proposition}\label{pn} 1) If $\det(\M)\ne 0$ then
the finite dimensional evolution algebra $\R^n$ has unique
absolute nilpotent $(0,...,0)$.

2) If $\det(\M)=0$ and rank$(\M)=n-1$ then the evolution algebra
$R^n$ has unique absolute nilpotent $(0,...,0)$ if and only if
\begin{equation}\label{con}
\det(\M_{i_0n})\cdot\det(\M_{n-1})>0,
\end{equation}
for some $i_0\in \{1,\dots,n-1\}$.
\end{proposition}
\proof 1) Straightforward.

2) If rank$(\M)=n-1$ then from (\ref{n11}) we get
\begin{equation}\label{n12}
x_i^2=-{\det(\M_{in})\over \det(M_{n-1})}x_n^2, \ \ i=1,\dots,n-1.
\end{equation}
From (\ref{n12}) it follows that the condition (\ref{con}) is
necessary and sufficient to have unique solution $(0,\dots,0)$.
\endproof

For a CEA $E^{[s,t]}$ with matrix $\M^{[s,t]}$ denote
$$\mathcal T_{nil}=\{(s,t)\in \mathcal T: E^{[s,t]}\ \ \mbox{has unique absolute
nilpotent}\}, \ \ \mathcal T_{nil}^0=\mathcal T\setminus \mathcal
T_{nil}.$$

The following theorem gives an answer on problem of existence of
``uniqueness of absolute nilpotent element'' property transition.

\begin{thm}\label{t5} 1) There are CEAs which have not ``uniqueness of
absolute nilpotent element'' property transition.

2) There is CEA which has ``uniqueness of absolute nilpotent
element'' property transition.
\end{thm}
\proof Denote  $ d(s,t)=\det(\M^{[s,t]})$. By equation (\ref{2}) we
get
\begin{equation}\label{dc}
d(s,t)=d(s,\tau)d(\tau,t), \ \ \mbox{for all} \ \ \tau, \ s<\tau<t.
\end{equation}
As it was mentioned above, the equation (\ref{dc}) is known as
Cantor's second equation.

1) The equation (\ref{dc}) has solutions $d(s,t)= {\Phi(s)\over
\Phi(t)}$, where $\Phi(t)\ne 0$ is an arbitrary function. Thus for
such solutions we conclude that if $d(s_0,t_0)\ne 0$ for some
$(s_0,t_0)$ then $d(s,t)\ne 0$ for any $(s,t)$. Consequently,
corresponding CEAs have not ``uniqueness of absolute nilpotent
element'' property transition.

2) Note that the equation (\ref{dc}) has solution $d(s,t)=f(t)$,
where $f(t)=1$ for $t<1$ and $f(t)=0$ otherwise. For this solution
we have $d(s,t)=1, s<t<1$ and $d(s,t)=0$, $t\geq 1$. For some
$t\geq 1$ one can construct a matrix $\M^{[s,t]}$ which does not
satisfy uniqueness condition mentioned in part 2) of Proposition
\ref{pn}. Indeed let us consider the matrix
$\M^{[s,t]}=\left(a^{[s,t]}_{ij}\right)_{i,j=1,2}$ with entries as
in (\ref{aa}). The second equation of the system (\ref{4}) has a
solution:
$$ \beta(s,t)=\left\{\begin{array}{ll}
1, \ \ \mbox{if} \ \ s<t<1\\[2mm]
0, \ \ \mbox{if} \ \ t\geq 1.\\[2mm]
\end{array}\right.
$$
Substituting this solution in the first equation of (\ref{4}) we
obtain
$$ \alpha(s,t)=\left\{\begin{array}{ll}
\psi(t)-\psi(s), \ \ \mbox{if} \ \ s<t<1\\[2mm]
g(t), \ \ \mbox{if} \ \ t\geq 1,\\[2mm]
\end{array}\right.
$$
where $\psi$ and $g$ are arbitrary functions. The corresponding
matrix has the following form
$$
\M^{[s,t]} ={1\over 2}\left(\begin{array}{cc}
2+\psi(t)-\psi(s)& -\psi(t)+\psi(s)\\[3mm]
\psi(t)-\psi(s) & 2-\psi(t)+\psi(s)
 \end{array}\right), \ \ \mbox{if} \ \ s<t<1,
$$ and
\begin{equation}\label{mm}
 \M^{[s,t]} ={1\over 2}\left(\begin{array}{cc}
1+g(t)& 1-g(t)\\[3mm]
1+g(t)& 1-g(t)
 \end{array}\right), \ \ \mbox{if} \ \ t\geq 1.
 \end{equation}
We have
$$ d(s,t)=\det(\M^{[s,t]})=\left\{\begin{array}{ll}
1, \ \ \mbox{if} \ \ s<t<1\\[2mm]
0, \ \ \mbox{if} \ \ t\geq 1.\\[2mm]
\end{array}\right.
$$
Assume $g(t)\ne -1$ then for (\ref{mm}) the equation (\ref{n12})
has the form
$$x_1^2=-{1-g(t)\over 1+g(t)}x_2^2.$$ This equation has infinitely
many solutions if $|g(t)|>1$ for some $t\geq 1$. Thus
corresponding CEA has ``uniqueness of absolute nilpotent element''
property transition.
\endproof

Now let us construct the set $\mathcal T_{nil}$ for Examples 1-5: It
is easy to see that
$$\det(\M^{[s,t]})=\begin{cases}
e^{-3At}, & \mbox{for Example 1};\\
(\lambda \mu)^{t}, & \mbox{for Example 2};\\
{\Phi(t)\over\Phi(s)}, & \mbox{for Example 3};\\
\prod_{i=1}^n \M^{[s,t]}_{ii} , & \mbox{for Example 4};\\
1,& \mbox{for Example 5}.
\end{cases}
$$
Thus in each one of the considered examples we have
$\det(\M^{[s,t]})\ne 0$, consequently, $\mathcal T_{nil}=\mathcal
T$, i.e. the CEAs constructed in Examples 1-5 have not
``uniqueness of the absolute nilpotent element'' property
transition.

There are CEAs which have infinitely many absolute nilpotent
elements independently on time. For example, take $\M^{[s,t]}$ with
identical rows\\ $({\Phi(s)\over \Phi(t)}, 0,0,\dots,0)$, where
$\Phi$ is an arbitrary function with $\Phi(t)\ne 0$ for all $t$. It
is easy to see that this matrix satisfies the equation (\ref{2}),
hence it determines a CEA, $E^{[s,t]}$, which has infinitely many
absolute nilpotent elements: $(0,x_2,\dots,x_n)$, where
$x_2,\dots,x_n\in \R$ are arbitrary numbers. Thus for this example
we have $\mathcal T_{nil}=\emptyset$, $\mathcal T_{nil}^0=\mathcal
T$. In other words the CEA has not ``non-uniqueness of absolute
nilpotent element'' property transition.

\begin{rk}
These examples (Examples 1-4) of ``uniqueness of nilpotent element''
property transition of CEAs with time-parameter are similar to the
``uniqueness of Gibbs phase'' property transition, i.e. phase
transition of physical systems with respect to
temperature-parameter, $T>0$. Usually there is a phase transition if
the temperature is very low ($T\sim 0$) or if it is very high
($T\sim +\infty$) (see \cite{ge}). Example 5 is an analogue of a
physical system which has unique (Gibbs) phase for any temperature.
There a lot of examples of such physical systems (see e.g.
\cite{ge}).
\end{rk}

\subsection{Idempotent elements transition}

A element $x$ of an algebra $\A$ is called {\it idempotent} if
$x^2=x$; such points of an evolution algebra are especially
important, because they are the fixed points (i.e. $V(x)=x$) of the
evolution operator $V$, (\ref{v}).  We denote by ${\mathcal Id}(E)$
the idempotent elements of an algebra $E$. Using (\ref{v}) the
equation $x^2=x$ can be written as
\begin{equation}\label{v1}
x_j= \sum_{i=1}^na_{ij}x_i^2, \ \ j=1,\dots,n.
\end{equation}
The general analysis of the solutions of the system (\ref{v1}) is
very difficult. We shall solve this problem for the CEA $E^{[t]}$,
$t\geq 0$, corresponding to the Example 2. In case of Example 2 the
system (\ref{v1}) has the following form
\begin{equation}\label{v2}\begin{cases}
2x=(\lambda^t+\mu^t)x^2+(\lambda^t-\mu^t)y^2;\\
2y=(\lambda^t-\mu^t)x^2+(\lambda^t+\mu^t)y^2,
\end{cases}
\end{equation}
where $\lambda>0, \ \ \mu>0$ and $t\geq 0$.

 {\it Case} $\lambda=\mu$. It is easy to see that if $\lambda=\mu$ then the
system (\ref{v2}) has only four solutions $0=(0,0),
z_1=z_1(t)=(0,\lambda^{-t}), z_2=z_2(t)=(\lambda^{-t},0),
z_3=z_3(t)=(\lambda^{-t},\lambda^{-t})$.

{\it Case} $\lambda\ne \mu$. For $\lambda\ne\mu$ the solutions $0$
and $z_3$ still exist. If $x=0$ or $y=0$ there is no any new
solution. Thus we consider the case $xy\ne0$. Denote
$$u={x\over y}, \ \ \gamma(t)={\lambda^t-\mu^t\over
\lambda^t+\mu^t}={(\lambda/\mu)^t-1\over (\lambda/\mu)^t+1}.$$ For
$t>0$ it is easy to see that if $\lambda<\mu$ then $-1<\gamma(t)<0$
and if $\lambda>\mu$ then $0<\gamma(t)<1$. Note that for $t=0$ there
is no any new solution. From system (\ref{v2}) we get

\begin{equation}\label{v3}
\gamma(t)u^3-u^2+u-\gamma(t)=(u-1)\left(\gamma(t)u^2+(\gamma(t)-1)u+\gamma(t)\right)=0.
\end{equation}

{\sl Subcase} $\lambda<\mu$. In this case for any $t>0$ the equation
(\ref{v3}) has three solutions
\begin{equation}\label{v4}
u_1=1, \ \
u_{\pm}={1-\gamma(t)\pm\sqrt{1-2\gamma(t)-3\gamma^2(t)}\over
2\gamma(t)}.
\end{equation}

{\sl Subcase} $\lambda>\mu$. In this case the number of solutions to
the equation (\ref{v3}) varies by $\gamma$, i.e.

\begin{equation}\label{v5}
\mbox{solutions to (\ref{v3})}=\begin{cases}
1,& \mbox{if} \ \ {1\over 3}\leq\gamma(t)<1;\\
1, u_-, u_+& \mbox{if} \ \ 0<\gamma(t)<{1\over 3},\\
\end{cases}
\end{equation}
where $u_{\pm}$ are defined in (\ref{v4}).

Now we shall describe $x, y$ corresponding to the solutions of
(\ref{v3}). The case $u=1$, i.e. $x=y$ does not give any new
solution. For $u=u_{\pm}$ we have $x=u_{\pm}y$, substituting this in
the second equation of (\ref{v2}) after simple calculations we get
the following two non-zero solutions to (\ref{v2}):
$$x_{\pm}={\mu^t\pm\sqrt{\lambda^t(2\mu^t-\lambda^t)}\over \mu^t(\lambda^t\pm\sqrt{\lambda^t(2\mu^t-\lambda^t)})},
\ \  y_{\pm}={\lambda^t-\mu^t\over
\mu^t(\lambda^t\pm\sqrt{\lambda^t(2\mu^t-\lambda^t)})}.$$ Note that
$x_{\pm}, \ \ y_{\pm}$ are well defined for any $\lambda\ne \mu$.
For $\lambda>\mu$ we have critical time
\begin{equation}\label{tt}
t_{\rm c}={\ln2\over \ln\lambda-\ln\mu},
\end{equation}
 which is the unique
solution to the equation $\gamma(t)={1\over 3}$.

Thus we have proved the following

\begin{proposition}\label{po} We have
$${\mathcal Id}(E^{[t]})=\begin{cases}
\{0, z_1,z_2,z_3\}, & \mbox{if} \ \ \lambda=\mu; \\
\{0,z_3, (x_-,y_-),(x_+,y_+)\}, & \mbox{if} \ \ \lambda<\mu; \\
\{0,z_3\}, & \mbox{if} \ \ \lambda>\mu; \ t\geq t_c \\
\{0,z_3, (x_-,y_-),(x_+,y_+)\}, & \mbox{if} \ \ \lambda>\mu, \ t< t_c. \\
\end{cases}$$
\end{proposition}

This proposition gives a very interesting ``a fixed set of
idempotent elements'' property transition, i.e. we have
\begin{cor} The CEA $E^{[t]}$ constructed in Example 2 has not ``a fixed set of idempotent
elements'' property transition if $\lambda \leq \mu$; it has such
property transition if $\lambda>\mu$. Moreover the transition point
(the critical time) is $t=t_c$ defined by formula (\ref{tt}).
\end{cor}

\begin{rk} There are exactly solvable models in statistical
mechanics, here an imprecise notion of ``exactly solvable'' as
meaning: ``The solutions can be expressed explicitly in terms of
some previously known functions'' is also sometimes used \cite{b}.
In such models, for example, the critical temperature can be
expressed explicitly. Comparing this with our examples of a property
transition we also can say that a property transition of a CEA is
exactly solvable if the critical time can be found exactly. Thus our
Example 2 is exactly solvable for investigation of properties of
idempotent elements.
\end{rk}

\section*{ Acknowledgements}

 The first and second authors were supported by Ministerio
de Ciencia e Innovaci\'on (European FEDER support included), grant
MTM2009-14464-C02, and by Xunta de Galicia, grant Incite09 207 215
PR. The third author thanks the Department of Algebra, University
of Santiago de Compostela, Spain,  for providing financial support
of his visit to the Department (February-April 2010). The authors
are grateful to both referees for helpful suggestions. The part 2)
of Theorem 4.14 and part 2) of Proposition \ref{pn} are due to a
suggestion of one referee.

{}

\begin{thebibliography}{99}
\bibitem{b} R.J. Baxter, \emph{Exactly solved models in statistical
mechanics}, Reprint of the 1982 original. Academic Press, Inc.
[Harcourt Brace Jovanovich, Publishers], London, (1989).

\bibitem{e1} I.M.H. Etherington, \emph{Genetic algebras}, Proc. Roy. Soc. Edinburgh. 59, 242--258 (1939).

\bibitem{e2} I.M.H. Etherington, \emph{Duplication of linear algebras}, Proc. Edinburgh Math. Soc. (2) 6, 222--230 (1941).

\bibitem{e3} I.M.H. Etherington, \emph{Non-associative algebra and the simbolism of genetics},
Proc. Roy. Soc. Edinburgh. 61, 24--42 (1941).

\bibitem{g} G. Galperin, A. Zemlyakov, \emph{Mathematical billiards}, Nauka, Moscow,
1990 (in Russian).

\bibitem{gp} N.N. Ganikhodjaev (Ganihodzhaev), \emph{On stochastic processes generated by quadric
operators},  J. Theoret. Probab.  4(4), 639--653 (1991).

\bibitem{gam} N.N. Ganikhodjaev,  H. Akin, F.M. Mukhamedov, \emph{On the
ergodic principle for Markov and quadratic stochastic processes and
its relations},  Linear Algebra Appl.  416(2-3), 730--741 (2006).

\bibitem{ge} H.-O. Georgii, \emph{Gibbs measures and phase transitions}, de
Gruyter Studies in Mathematics, 9. Walter de Gruyter  \& Co.,
Berlin, (1988).

\bibitem{ht} P. H\"anggi, H. Thomas, \emph{Time evolution, correlations,
and linear response of non-Markov processes}, Zeitschrift Phys. B.
26, 85--92 (1977).

\bibitem{k} A.N. Kolmogorov, \emph{On analytic methods in probability
theory}, Uspekhi Mat. Nauk. 5, 5--41 (1938) (Russian), English
transl., Selected works of A.N. Kolmogorov, Vol. II, Kluwer,
Dordrecht 1992, article 9.

\bibitem{ly} Y.I. Lyubich, \emph{Mathematical structures in population
genetics}, Springer-Verlag, Berlin, 1992.

\bibitem{m} M.L. Reed, \emph{Algebraic structure of genetic inheritance},  Bull. Amer. Math.
Soc. (N.S.)  34(2), 107--130  (1997).

\bibitem{rt} U.A. Rozikov, J. P. Tian, \emph{Evolution algebras generated by Gibbs measures}, Preprint ICTP,
IC2009013.

\bibitem{gs} T.A. Sarymsakov, N.N. Ganikhodjaev, \emph{Analytic methods in the
theory of quadric stochastic operators},  J. Theoret. Probab.  3(1),
51--70 (1990).

\bibitem{t} J. P. Tian, \emph{Evolution algebras and their applications},
Lecture Notes in Mathematics, 1921, Springer-Verlag, Berlin, 2008.

\bibitem{w} A. W\"orz-Busekros, \emph{Algebras in genetics}, Lecture Notes in
Biomathematics, 36. Springer-Verlag, Berlin-New York, 1980.
\end{thebibliography}
\end{document}